%May 8, 2004
\input amstex
\documentstyle{amsppt}
\magnification=1200
\hoffset=-0.5pc
\nologo
\vsize=57.2truepc
\hsize=38.5truepc

\spaceskip=.5em plus.25em minus.20em
\define\din{d}

\define\Bobb{\Bbb}
\define\fra{\frak}
\define\HG{G}
\define\hg{g}
\define\GHH{H}
\define\KK{H}
\define\KKK{H}
\define\kk{h}

\define\GH{H}

\define\sk{k}
\define\armcusgo{1}
\define\atiybern{2}
\define\chapuone{3}
\define\farakora{4}
\define\freudone{5}
\define\hartsboo{6}
\define\howeone{7}
\define\poiscoho{8}
\define\souriau{9}
\define\oberwork{10}
\define\kaehler{11}
\define\kaehredu{12}
\define\lradq{13}
\define\jacobstw{14}
\define\kostasix{15}
\define\lanmaone{16}
\define\lermonsj{17}
\define\mccrione{18}
\define\mezzthom{19}
\define\satakboo{20}
\define\scorzone{21}
\define\scorztwo{22}
\define\sekiguch{23}
\define\severone{24}
\define\zakone{25}
\define\zaktwo{26}
\topmatter
\title Singular Poisson-K\"ahler geometry of Scorza
varieties and their secant varieties
\endtitle
\author Johannes Huebschmann
\endauthor
\affil
Universit\'e des Sciences et Technologies
de Lille
\\
U. F. R. de Math\'ematiques
\\
CNRS-UMR 8524
\\
F-59 655 VILLENEUVE D'ASCQ, France
\\
Johannes.Huebschmann\@math.univ-lille1.fr
\endaffil
\date{May 8, 2004}
\enddate
\abstract{Each Scorza variety and its secant varieties in the
ambient projective space are identified, in the realm of singular
Poisson-K\"ahler geometry, in terms of projectivizations of
holomorphic nilpotent orbits in suitable Lie algebras of hermitian
type, the holomorphic nilpotent orbits, in turn, being affine
varieties. The ambient projective space acquires an exotic
K\"ahler structure, the closed stratum being the Scorza variety
and the closures of the higher strata its secant varieties. In
this fashion, the secant varieties become exotic projective
varieties. In the rank 3 case, the four regular Scorza varieties
coincide with the four critical Severi varieties. In the standard
cases, the Scorza varieties and their secant varieties arise also
via K\"ahler reduction. An interpretation in terms of constrained
mechanical systems is included.}
\endabstract

\address{\smallskip
\noindent
USTL, UFR de Math\'ematiques, CNRS-UMR 8524
\newline\noindent
F-59 655 Villeneuve d'Ascq C\'edex,
France
\newline\noindent
Johannes.Huebschmann\@math.univ-lille1.fr}
\endaddress
\endkeywords

\endtopmatter
\document
\leftheadtext{Johannes Huebschmann} \rightheadtext{Singular
geometry of Scorza varieties and their secant varieties}

\beginsection 1. Introduction

Let $m \geq 2$. A {\it Severi variety\/} is a non-singular variety
$X$ in complex projective $m$-space $\Bobb P^m \Bobb C$ having the
property that, for some point $O \not \in X$, the projection from
$X$ to $\Bobb P^{m-1} \Bobb C$ is a closed immersion, cf.
\cite\hartsboo\ (Ex.~3.11 p. 316),\,
\,\cite\lanmaone,\,\cite\zakone. Let $X$ be a Severi variety and
$n$ the dimension of $X$. The critical cases are when $m = \frac
32 n + 2$. Zak \cite\zakone\ proved that only the following {\it
four\/} critical cases occur:
\newline\noindent
(1.1) $X=\Bobb P^2 \Bobb C \subseteq \Bobb P^5 \Bobb C$ (Veronese
embedding)
\newline\noindent
(1.2) $X=\Bobb P^2 \Bobb C \times \Bobb P^2 \Bobb C \subseteq
\Bobb P^8 \Bobb C$ (Segre embedding)
\newline\noindent
(1.3) $X=\roman G_2(\Bobb C^6) = \roman U(6)\big /(\roman
U(2)\times \roman U(4)) \subseteq \Bobb P^{14} \Bobb C$ (Pl\"ucker
embedding)
\newline\noindent
(1.4) $X=\roman{Ad}(\fra e_{6(-78)}) \big/(\roman{SO}(10,\Bobb R)
\cdot \roman{SO}(2,\Bobb R)) \subseteq \Bobb P^{26} \Bobb C$.
\newline\noindent
These varieties arise from the projective planes
over the {\it four\/} real normed division algebras
(reals, complex numbers, quaternions, octonions)
by complexification, cf. e.~g. \cite\atiybern.

In the book \cite\zaktwo, Zak introduced and classified {\it
Scorza varieties\/}. These generalize the critical Severi
varieties. Recall that, given a non-singular variety $X$ in $\Bobb
P^m \Bobb C$, for $0 \leq k \leq m$, the $k$'th {\it secant
variety\/} $S^k(X)$ is the projective variety in $\Bobb P^m \Bobb
C$ which arises as the closure of the union of all $k$-dimensional
projective spaces in $\Bobb P^m \Bobb C$ that contain $k+1$
independent points of $X$. Then $S^0(X) = X$ and $S^1(X)$ is the
ordinary secant variety, referred to as well as {\it chordal
variety\/}. For $k \geq 2$, a $k$-{\it Scorza variety\/} is
defined to be a non-singular complex projective variety of maximal
dimension among those varieties whose $(k-1)$-secant variety is
not the entire ambient projective space; a precise description of
that notion of maximality will be recalled in Section 4 below.

The Severi variety (1.4) is a very special 2-Scorza variety. Let
$\sk \geq 2$. We now list the other $\sk$-Scorza varieties, cf.
Theorem 5.6 in Chap. 6 of \cite\zaktwo:

\noindent (1.1.$\sk$) $X=\Bobb P^{\sk} \Bobb C \subseteq \Bobb
P^{\frac{\sk(\sk+3)}2} \Bobb C$ (Veronese embedding)
\newline\noindent
(1.2.$\sk$.r) $X=\Bobb P^{\sk} \Bobb C \times \Bobb P^{\sk} \Bobb
C \subseteq \Bobb P^{\sk(\sk+2)} \Bobb C$ (Segre embedding)
\newline\noindent
(1.2.$\sk$.n) $X=\Bobb P^{\sk} \Bobb C \times \Bobb P^{\sk+1}
\Bobb C \subseteq \Bobb P^{\sk^2+3\sk+1} \Bobb C$ (Segre
embedding)
\newline\noindent
(1.3.$\sk$.r) $X=\roman G_2(\Bobb C^{2(\sk+1)}) \subseteq \Bobb
P^{\sk(2\sk+3)} \Bobb C$ (Pl\"ucker embedding)
\newline\noindent
(1.3.$\sk$.n) $X=\roman G_2(\Bobb C^{2\sk+3)}) \subseteq \Bobb
P^{2\sk^2+5\sk + 2} \Bobb C$ (Pl\"ucker embedding)
\newline
For reasons that will become clear below we will refer to the
Scorza varieties (1.1.$\sk$), (1.2.$\sk$.r), (1.3.$\sk$.r), and
(1.4) as {\it regular\/} and to the remaining ones as {\it
non-regular\/}. The critical Severi varieties are exactly the
regular 2-Scorza varieties.

When the chordal variety of a non-singular projective variety $Q$
in $\Bobb P^5 \Bobb C$ is a hypersurface (and not the entire
ambient space) the projection from a generic point gives an
embedding in $\Bobb P^4 \Bobb C$. A classical result of Severi
\cite\severone\ says that the Veronese surface is the only surface
(not contained in a hyperplane) in $\Bobb P^5 \Bobb C$ with this
property. This is the origin of the terminology \lq\lq Severi
variety\rq\rq. In \cite\atiybern, the chordal varieties are
written as $Z_n(C)$ ($n=0,1,2,3$). The terminology \lq\lq Scorza
variety\rq\rq\ has been introduced in \cite\zaktwo, to honor
Scorza's pioneering work on linear normalization of varieties of
small codimension \cite\scorzone,\, \cite\scorztwo.

The purpose of the present paper is to exhibit an interesting
geometric feature of Scorza varieties, in particular of the
critical Severi varieties, in the world of singular
Poisson-K\"ahler geometry, to be given as Theorem 1.5 below: There
are exactly {\it four\/} simple regular rank 3 hermitian Lie
algebras over the reals; these result from the euclidean Jordan
algebras of hermitian ($3\times 3$)-matrices over the {\it four\/}
real normed division algebras by the superstructure construction;
and {\it each critical Severi variety arises from the minimal
holomorphic nilpotent orbit in such a Lie algebra.\/} The
resulting geometric insight into the four critical Severi
varieties will be made explicit in the Addendum 1.7 below. For
higher rank, that is, when the rank $r$ (say) is at least equal to
$4$, there are exactly {\it three\/} simple regular rank $r$
hermitian Lie algebras over the reals; these result from the
euclidean Jordan algebras of hermitian ($r\times r$)-matrices over
the {\it three\/} {\it associative\/} real normed division
algebras by the superstructure construction. Another result of
Zak's \cite\zaktwo, combined with Theorem 3.3.11 in \cite\kaehler,
entails that {\it each $\sk$-Scorza variety arises from the
minimal holomorphic nilpotent orbit in a simple hermitian Lie
algebra of rank $\sk+1$\/} $(\sk \geq 2)$; furthermore, the {\it
regular Scorza varieties arise in this fashion from regular simple
hermitian Lie algebas and hence, cf. Remark {\rm 5.10} and Theorem
{\rm 5.11} in Chap. 6 of \cite\zaktwo, from simple complex Jordan
algebras, and the non-regular Scorza varieties arise from
non-regular simple hermitian Lie algebas and hence from simple
positive definite hermitian Jordan triple systems\/}. In
particular, every simple complex Jordan algebra of rank at least 3
occurs here, indeed, the classification of regular Scorza
varieties parallels that of simple complex Jordan algebras and
hence that of tube domains. However not every positive definite
hermitian Jordan triple system gives rise to a Scorza variety; see
Remark 5.7 below for details.

We will spell out the resulting geometric insight into general
Scorza varieties, in the realm of singular Poisson-K\"ahler
geometry, in Theorem 1.5 below. To explain this geometric insight,
we need some preparation. Let $Y$ be a stratified space. A complex
analytic stratified K\"ahler structure on $Y$ in the sense of
\cite\kaehler\ is a stratified symplectic structure together with
a compatible complex analytic structure which, on each stratum,
combine to a K\"ahler structure; when at least two strata are
involved, we refer to an {\it exotic K\"ahler structure\/}. A
complete definition  will be reproduced in Section 5 below.
Suffice it to mention at this stage that the structure includes a
{\it real\/} Poisson algebra of continuous functions on $Y$ which,
on each stratum, restricts to an ordinary smooth symplectic
Poisson algebra. This Poisson structure is {\it independent\/} of
the complex analytic structure.

\proclaim{Theorem 1.5} Let $\sk \geq 2$, let $X$ be a
$\sk$-Scorza variety, and let $\Bobb P^m \Bobb C$ be the ambient complex
projective space. Then this projective space carries an exotic normal
K\"ahler structure with the following properties:
\newline\noindent
{\rm (1)} The closures of the strata constitute an ascending
sequence
$$
Q_1 \subseteq Q_2 \subseteq  \ldots \subseteq Q_{\sk+1}=\Bobb P^m \Bobb
C \tag1.6
$$
of normal K\"ahler spaces where the closed stratum $Q_1$ coincides
with the given Scorza variety and where, complex algebraically,
for $1 \leq \rho \leq \sk$, $Q_\rho$ is a projective determinantal
variety in $Q_{\sk+1} = \Bobb P^m \Bobb C$; in particular, in the
regular case, $Q_{\sk}$ is a projective degree $\sk+1$
hypersurface.
\newline\noindent
{\rm (2)}  For $1 \leq \rho \leq {\sk}$, $Q_{\rho+1}$ is the
$\rho$'th secant variety
$S^\rho(Q_1)$ of $Q_1$ in $Q_{\sk+1} = \Bobb P^m \Bobb C$.
\newline\noindent
{\rm (3)} For $2 \leq \rho \leq {\sk+1}$, $Q_{\rho-1}$ is the singular locus
of $Q_\rho$, in the sense of stratified K\"ahler spaces.
\newline\noindent
{\rm (4)} The exotic K\"ahler structure on $\Bobb P^m \Bobb C$
restricts to an ordinary K\"ahler structure on  $Q_1$ inducing,
perhaps up to rescaling, the standard hermitian symmetric space
structure.
\endproclaim

We now spell out the result explicitly in the regular rank 3 case,
which is somewhat special. This case corresponds to the critical
Severi varieties and includes the exceptional case of the Severi
variety (1.4) arising from the octonions (see below).

\proclaim{Addendum 1.7} For $m=5,8,14,26$, the complex projective
space $\Bobb P^m \Bobb C$ carries an exotic normal K\"ahler
structure with the following properties:
\newline\noindent
{\rm (1)}
The closures of the strata constitute an ascending sequence
$$
Q_1 \subseteq Q_2 \subseteq  Q_3=\Bobb P^m \Bobb C \tag1.6.SEVERI
$$
of normal K\"ahler spaces where, complex algebraically,
$Q_1$ is a (critical) Severi variety and
$Q_2$ a projective cubic hypersurface,
the chordal variety of $Q_1$.
\newline\noindent
{\rm (2)}
The singular locus of $Q_3$, in the sense of stratified
K\"ahler spaces, is the hypersurface $Q_2$, and that of $Q_2$
(still in the sense of stratified K\"ahler spaces)
is the non-singular variety $Q_1$; furthermore,  $Q_1$ is as well the
complex algebraic singular locus of $Q_2$.
\newline\noindent
{\rm (3)}
The exotic  K\"ahler structure
on $\Bobb P^m \Bobb C$ restricts to an ordinary K\"ahler structure
on  $Q_1$ inducing, perhaps up to rescaling,
the standard hermitian symmetric space structure.
\endproclaim

In case (1.4), the term \lq\lq determinantal variety\rq\rq\ refers
to Freudenthal's determinant, see Sections 8 and 10 of
\cite\kaehler.

For a compact K\"ahler manifold $N$ which is complex analytically
a projective variety, with reference to the Fubini-Study metric,
the Kodaira embedding will not in general be symplectic; the above
theorem shows that there are interesting situations where the {\sl
ambient complex  projective space carries an exotic K\"ahler
structure\/} which, via the {\sl Kodaira embedding, restricts to
the K\"ahler structure on\/} $N$. More generally, for a projective
variety $N$ with singularities, the correct question to ask is
whether, via the Kodaira embedding, $N$ inherits, from a suitable
exotic K\"ahler structure on the ambient complex projective space,
a (stratified) K\"ahler structure. In  \cite\kaehredu\ we have
developped a K\"ahler quantization scheme for not necessarily
smooth stratified K\"ahler spaces, including examples of K\"ahler
quantization on projective varieties with stratified K\"ahler
structure. This procedure applies to the circumstances of the
above theorem. It is worthwhile noting that ignoring the lower
strata means working on a non-compact space, and this then leads
to inconsistencies in the sense that the principle that
quantization commutes with reduction is violated. See
\cite\kaehredu\ (4.12) for details.

An interesting feature of the situations isolated in Theorem 1.5
is that the real Poisson algebra on $\Bobb P^m \Bobb C$ (beware:
it contains more functions than just the ordinary smooth ones)
detects the lower strata $Q_s \setminus Q_{s-1}$ ($1 \leq s \leq
\sk+1$ where $Q_0$ is understood to be empty) by means of the rank
of the Poisson structure, independently of the complex analytic
structures. This notion of rank can be made precise by means of
the Lie-Rinehart structure related with a general not necessarily
smooth Poisson structure which we introduced in \cite\poiscoho.

More geometric consequences will be explained later. In
particular, in Section 6 below we will show how the ascending
sequence (1.6) (and hence the ascending sequence (1.6.SEVERI) for
the classical cases (1.1), (1.2), (1.3)) results from K\"ahler
reduction which, in particular, exhibits, for $1 \leq s \leq
\sk+1$, each constituent $Q_s$ complex analytically as a
G(eometric)I(nvariant)T(heory)quotient. For $s=1$, this
construction includes Zak's Veronese mapping \cite\zaktwo, see
Remark 6.8 below. When the GIT-quotient is related with the
corresponding symplectic quotient via the nowadays familiar
Kempf-Ness observation, Theorem 1.5 provides geometric insight
into the singular structure of the corresponding symplectic
quotients. As for the exceptional case (1.4), we do not know
whether the ascending sequence (1.6.SEVERI) may be obtained from a
smooth K\"ahler manifold via K\"ahler reduction.

The underlying real spaces for the affine situation from which the
sequence (1.6) in case (1.1.$\sk$) arises by projectivization
result from symplectic reduction of the phase space of $n=\sk+1$
particles moving in $\Bobb R^n$, with reference to total angular
momentum zero. In \cite\lermonsj\ it has been shown that, in this
case, the $\roman{Sp}(n,\Bobb R)$-momentum mapping identifies the
reduced space, as a stratified symplectic space, with the closure
of a certain nilpotent orbit in $\fra{sp}(n,\Bobb R)$ and,
likewise, the lower strata are similar nilpotent orbits which
correspond to the reduced systems of $n$ particles moving in
$\Bobb R^k$ for $k<n$. The requisite stratified symplectic Poisson
algebra on the reduced spaces is that constructed in
\cite\armcusgo. Comparing the result in \cite\lermonsj\ just
quoted with the Kempf-Ness observation in geometric invariant
theory I noticed that the nilpotent orbits isolated in
\cite\lermonsj\ are precisely those, for the special case of
$\fra{sp}(n,\Bobb R)$ ($n \geq 1$), which I identified as
holomorphic nilpotent orbits in \cite\kaehler. This was the
starting point of the entire research program the present paper is
part of.

For a fixed rank $\sk+1$, in the two other (classical) cases, the
constituents of the ascending sequences (1.6) admit presumably a
similar interpretation in terms of certain  constrained systems in
mechanics;  see Remark 6.7 below for details. Here the complex
analytic structure does not seem to have a direct mechanical
significance; it helps understanding the kinematical description.
For issues related with quantization, the complex analytic
structure has its intrinsic significance, though, cf.
\cite\kaehredu.

I am indebted to F. Hirzebruch for having introduced me into
Severi varieties, to a competent referee who suggested that the
approach, spelled out in a draft of the manuscript only for Severi
varieties, should extend to general Scorza varieties, and to
J.~P.~ Chaput and F.~L.~Zak for some illuminating e-mail
correspondence about Scorza varieties.

\medskip\noindent {\bf 2. Lie algebras of hermitian type}
\smallskip\noindent
Following \cite\satakboo\ (p.~54), we define a (semisimple) Lie algebra
of {\it hermitian type\/} to be a pair $(\fra g, z)$ which consists of a
real semisimple Lie algebra $\fra g$ with a Cartan decomposition
$\fra g = \fra k \oplus \fra p$ and a central element $z$ of $\fra k$,
referred to as an $H$-{\it element\/}, such that
$J_z = \roman{ad}(z)\big |_{\fra p}$ is a (necessarily $K$-invariant)
complex structure on $\fra p$. Slightly more generally,
a {\it reductive Lie algebra of hermitian type\/} is a reductive Lie algebra
$\fra g$ together with an element $z \in \fra g$
whose constituent $z'$ (say) in the semisimple part
$[\fra g,\fra g]$ of $\fra g$ is an $H$-element for $[\fra g,\fra g]$
\cite\satakboo\ (p.~92). Below we will sometimes refer to $\fra g$ alone
(without explicit choice of $H$-element $z$) as a {\it hermitian\/} Lie
algebra. For a real semisimple Lie algebra $\fra g$, with Cartan
decomposition $\fra g = \fra k \oplus \fra p$, we write $G$ for an
appropriate Lie group having $\fra g = \roman {Lie}(G)$
(matrix realization or adjoint realization; both will do)
and $K$ for the (compact) connected subgroup of $G$ with \
$\roman{Lie}(K) = \fra k$; the requirement that $(\fra g,z)$ be of
hermitian type is equivalent to $G\big / K$ being a (non-compact)
hermitian symmetric space with complex structure induced by $z$.

A real semisimple hermitian Lie algebra $\fra g$ decomposes as
$\fra g = \fra g_0 \oplus \fra g_1 \oplus \ldots  \oplus\fra g_k$
where $\fra g_0$ is the maximal compact semisimple ideal and where
$\fra g_1,\dots,  \fra g_k$ are non-compact and simple.
For a non-compact simple Lie algebra with Cartan decomposition
$\fra g = \fra k \oplus \fra p$, the  $\fra k$-action on $\fra p$
coming from the adjoint representation of $\fra g$
is faithful and irreducible whence the center of $\fra k$
then is at most one-dimensional; indeed $\fra g$ has an $H$-element
turning it into a Lie algebra of hermitian type if and only if
the center of $\fra k$ has dimension one. In view of E.~Cartan's
infinitesimal classification of irreducible hermitian symmetric spaces,
the Lie algebras $\fra{su}(p,q)$ ($A_n, n \geq 1$, where $n+1 = p+q$),
$\fra{so}(2,2n-1)$ ($B_n, n \geq 2 $),
$\fra{sp}(n,\Bobb R)$ ($C_n, n \geq 2$),
$\fra{so}(2,2n-2)$ ($D_{n,1}, n>2$),
$\fra{so}^*(2n)$ ($D_{n,2}, n >2$)
together with the real forms
$\fra{e}_{6(-14)}$ and $\fra{e}_{7(-25)}$
of $\fra{e}_{6}$ and $\fra{e}_{7}$, respectively,
constitute a complete list of simple hermitian Lie algebras.

We refer to $(\fra g,z)$ as {\it regular\/} when the relative root
system is of type $\roman C_r$, $r \geq 1$; see \cite\kaehler\
(Proposition 3.3.2) for details. The natural number $r$ is then
the {\it real\/} rank of $\fra g$. Thus $\fra {sp}(1,\Bobb R)
\cong\fra {su}(1,1)\cong \fra {so}(2,1)$ and $\fra{so}^*(4)$ are
the only regular rank 1 simple hermitian Lie algebras; the $\fra
{so}(p,2)$'s ($p \geq 3$) are the only regular rank 2 simple
hermitian Lie algebras; in particular, $\fra {sp}(2,\Bobb R)\cong
\fra {so}(2,3)$, $\fra {su}(2,2)\cong \fra {so}(2,4)$,
$\fra{so}^*(8)\cong \fra {so}(2,6)$; $\fra {sp}(3,\Bobb R)$, $\fra
{su}(3,3)$, $\fra{so}^*(12)$ and $\fra{e}_{7(-25)}$ are the only
regular rank 3 simple hermitian Lie algebras; and, for $n \geq 4$,
$\fra {sp}(n,\Bobb R)$, $\fra {su}(n,n)$, and $\fra{so}^*(4n)$ are
the only regular rank $n$ simple hermitian Lie algebras. Moreover,
$\fra{so}^*(10)$, $\fra{e}_{6(-14)}$, and the $\fra {su}(p,2)$'s
for $p>2$ are the only non-regular rank 2 simple hermitian Lie
algebras; and, for $n \geq 3$, the $\fra {su}(p,n)$'s for $p>n$
and $\fra{so}^*(4n+2)$ constitute a complete list of the
non-regular rank $n$ simple hermitian Lie algebras.

\medskip\noindent{\bf 3. Jordan algebras and hermitian Jordan triple systems}
\smallskip\noindent
The symmetric constituent $\fra p$ of the Cartan decomposition
$\fra g =\fra k \oplus \fra p$ of a Lie algebra $(\fra g,z)$ of
hermitian type inherits a hermitian Jordan triple system structure
which, in case $\fra g$ is simple and regular, comes down to a
simple complex Jordan algebra. We now explain briefly the relevant
pieces of structure which we shall need; see Section 7 of
\cite\kaehler\ for more details and notation.

As usual, denote by $\Bobb O$ the octonions or Cayley numbers. For
$\Bobb K = \Bobb R,\,\Bobb C,\,\Bobb H,\,\Bobb O$ and $n \geq 1$,
let $\Cal H_n(\Bobb K)$ be the {\it euclidean Jordan algebra\/} of
hermitian $(n \times n)$-matrices over $\Bobb K$ ($\Cal H_n(\Bobb
R) = \roman S^2_{\Bobb R}[\Bobb R^n]$), with Jordan product
$\,\circ\,$ given by $x \circ y = \frac 12 (xy + yx)$ ($x,y \in
\Cal H_n(\Bobb K)$) where $n \leq 3$ when $\Bobb K = \Bobb O$. See
\cite{\farakora,\,\mccrione,\,\satakboo} for notation and details.
Here $\Cal H_1(\Bobb R)\cong\Cal H_1(\Bobb C) \cong\Cal H_1(\Bobb
H)\cong\Cal H_1(\Bobb O)\cong \Bobb R$, $\Cal H_2(\Bobb O)$ is
isomorphic to the euclidean Jordan algebra $J(1,9)$ arising from
the Lorentz form of type $(1,9)$ on $\Bobb R^{10}$, and $\Cal
H_3(\Bobb O)$ is the real exceptional rank 3 Jordan algebra of
dimension 27. Complexification yields the complex simple Jordan
algebras $\roman S_{\Bobb C}^2[\Bobb C^n]$, $\roman M_{n,n}(\Bobb
C)$, $\Lambda^2[\Bobb C^{2n}]$, $\Cal H_3(\Bobb O_{\Bobb C})=\Cal
H_3(\Bobb O)\otimes \Bobb C$, and it is well known that, when $n
\geq 3$, these exhaust the simple complex Jordan algebras of rank
$ \geq 3$; see e.~g. \cite\farakora.

Let $(\fra g,z)$ be a simple Lie algebra  of hermitian type of
rank $r\geq 3$. The Cartan decomposition $\fra g = \fra k \oplus
\fra p$ has the following form where the decomposition (C.2) is
spelled out for the reductive hermitian Lie algebra $\fra u(q,p)$
instead of its simple brother $\fra {su}(q,p)$:
\newline\noindent
(C.1) $\fra{sp}(n,\Bobb R)= \fra u(n) \oplus \roman S_{\Bobb
C}^2[\Bobb C^n]$ ($r=n\geq 3$)
\newline\noindent
(C.2) $\fra u(q,p)=(\fra u(q)\oplus\fra u(p)) \oplus \roman
M_{q,p}(\Bobb C)$ ($p\geq q=r\geq 3$)
\newline\noindent
(C.3) $\fra {so}^*(2n)= \fra u(n) \oplus \Lambda^2[\Bobb C^{n}]$
($n\geq 6$, $r=\left[\frac n2\right]$)
\newline\noindent
(C.4) $\fra e_{7(-25)} =(\fra e_{6(-78)}\oplus\Bobb R)\oplus\Cal
H_3(\Bobb O_{\Bobb C})$, $\fra e_{6(-78)}$ being the compact form
of $\fra e_6$ (only for $r=3$).
\newline\noindent
Let $K$ be the compact constituent in the Cartan decomposition of
the adjoint group $\HG$ of $\fra g$, and let $K^{\Bobb C}$ be its
complexification.
 The resulting
(unitary) $K$-representation on the complex vector space $\fra p$
extends to a $K^{\Bobb C}$-representation on $\fra p$ which turns
the latter into a pre-homogeneous space which we refer to as {\it
regular\/} whenever the corresponding Lie algebra of hermitian
type is regular. Explicitly, these representations have the
following forms:
\newline\noindent
(R.1) The symmetric square of the standard $\roman {GL}(n,\Bobb
C)$-representation on $\Bobb C^n$, given by $x\cdot S = x S x^t$,
for $x \in \roman{GL}(n,\Bobb C)$ and $S \in \roman S_{\Bobb
C}^2[\Bobb C^n]$; it is always regular.
\newline\noindent
(R.2) The standard representation of $\roman {GL}(q,\Bobb C)
\times \roman {GL}(p,\Bobb C)$ on $\roman M_{q,p}(\Bobb C)$ given
by
$$
(x,y)\cdot M = x M y^t,\quad x \in \roman{GL}(p,\Bobb C), \ y \in
\roman{GL}(q,\Bobb C), \ M \in \roman M_{q,p}(\Bobb C);
$$
it is regular if and only if $p=q$.

\noindent (R.3) The exterior square of the standard $\roman
{GL}(n,\Bobb C)$-representation on $\Bobb C^{n}$; it is regular if
and only if $n$ is even.
\newline\noindent
(R.4) The classical representation of $\roman E_6(\Bobb C)$,
extended by a central copy of $\Bobb C^*$; it is regular and has
complex dimension 27. This representation was studied already by
E. Cartan.

The  complexification $\fra g^{\Bobb C} =  \fra p^+\oplus \fra
k^{\Bobb C} \oplus \fra p^-$, where $\fra p^+$ and $\fra p^-$
refer to the holomorphic and antiholomorphic constituents of $\fra
p ^{\Bobb C}$, respectively, inherits a real symmetric Lie algebra
structure, the requisite involution being complex conjugation, and
this structure, in turn, determines the structure of a simple
positive definite hermitian Jordan triple system (JTS) on $\fra
p^+$. Indeed, every positive definite hermitian JTS arises in this
fashion, cf. \cite\satakboo\ (Proposition II.3.3 on p.~56). Under
this correspondence, the regular simple Lie algebras of hermitian
type correspond to positive definite hermitian JTS's which arise
from a Jordan algebra (with unit element). Furthermore, the action
on $\fra p^+$, of the compact constituent $K$ in the Cartan
decomposition of the adjoint group $\HG$ of $\fra g$, extends to
an action of the complexification $K^{\Bobb C}$ on $\fra p^+$ and
turns the latter into a pre-homogeneous space for $K^{\Bobb C}$.
The obvious isomorphism of complex vector spaces between $\fra p$
and $\fra p^+$ identifies the pre-homogeneous space structures. In
this fashion, $\fra p$ acquires a JTS structure. More details may
be found in Sections 7 and 8 of \cite\kaehler.

\medskip\noindent
{\bf 4. Jordan rank and Severi, Scorza, and secant varieties}
\smallskip\noindent
Let $\fra g$ still be a simple hermitian Lie algebra, with Cartan
decomposition $\fra g = \fra k \oplus \fra p$, and let $r$ be the
real rank of $\fra g$. For $1 \leq s \leq r$, let $\Cal O_s
\subseteq \fra p$ be the subspace of Jordan rank $s$. This notion
of Jordan rank is discussed in Section 7 of \cite\kaehler. In the
regular (i.~e. Jordan algebra) case, the Jordan rank of a non-zero
element is the number of non-zero eigenvalues in its spectral
decomposition, with multiplicities counted \cite\farakora\
(p.~77). For $\fra g = \fra{sp}(r,\Bobb R)$ and $\fra g =
\fra{su}(p,r)$ ($p\geq r \geq 1$), the Jordan rank amounts to the
ordinary rank of a matrix, where $\roman S_{\Bobb C}^2[\Bobb C^r]$
is identified with the symmetric complex  $(r\times r)$-matrices.
For $\fra g = \fra{so}^*(2n)$ ($r=\left[\frac n2\right]\geq 2$),
when we identify $\Lambda^2[\Bobb C^{n}]$ with the skew-symmetric
complex $(n\times n)$-matrices, the Jordan rank amounts to one
half the ordinary rank of a matrix. The resulting decomposition of
$\fra p$ is a stratification whose strata coincide with the
$K^{\Bobb C}$-orbits, and the closures constitute an ascending
sequence
$$
\{0\} \subseteq \overline {\Cal O_1} \subseteq \overline {\Cal
O_2} \subseteq \ldots \subseteq \overline {\Cal O_r} = \fra p
\tag4.1
$$
of complex affine varieties. Projectivization yields the ascending sequence
$$
Q_1 \subseteq Q_2 \subseteq \ldots \subset Q_r =\Bobb P(\fra p).
\tag4.2
$$
As far as the complex analytic structures are concerned, this is
the sequence (1.6) when $r \geq 3$, in particular the sequence
(1.6.SEVERI) for $r=3$ in the regular case.

In the Severi case, that is, under the circumstances of Theorem
1.5 for for $k=2$ in the regular case or, equivalently, under the
circumstances of the Addendum 1.7, by construction, $Q_2$ is a
determinantal cubic, the requisite determinant over the octonions
being that introduced by Freudenthal \cite\freudone, and $Q_1$ is
the corresponding Severi variety. In fact, $Q_1$ is the {\it
closed\/} $K^{\Bobb C}$-orbit in $\Bobb P(\fra p)$, and the
homogeneous space descriptions (1.1)--(1.4) are immediate. For
$\fra g =\fra e_{7(-25)}$, the cubic $Q_2$ is the (projective)
generic norm hypersurface, sometimes referred to in the literature
as {\it Freudenthal\/} cubic; it has been studied  already by E.
Cartan, though. Jacobson has shown that this cubic is rational
\cite\jacobstw.

We now recall a description of  Scorza varieties. Let $X\subseteq
\Bobb P^m \Bobb C$ be a non-singular complex projective variety
which does not lie in a hyperplane. Let $k_0(X)$ be the smallest
natural number $k$ such that $S^k(X)=\Bobb P^m \Bobb C$. Given
another projective variety $Y$, let $\delta(X,Y) =\dim X
+\dim(Y)+1-\dim S(X,Y)$ be the {\it defect\/} of $X$ and $Y$ where
$S(X,Y)$ refers to the join of $X$ and $Y$, that is, to the
closure of the union of the lines through a point of $X$ and a
point of $Y$; following Zak \cite\zaktwo, for $1 \leq i \leq
k_0(X)$, let
$$
\delta_i=\delta_i(X) = \delta(X,S^{i-1}(X))
= \dim X + \dim S^{i-1}(X) +1 -  \dim S^i(X),
$$
and write $\delta = \delta_1$. Zak proved that $\delta_i \geq
\delta_{i-1} + \delta$ whenever $2 \leq i \leq k_0(X)$
(\cite\zaktwo, Chap. 5, Theorem 1.8) and deduced that $k_0(X) \leq
\frac{\dim X} \delta$. Rather than reproducing Zak's definition of
a Scorza variety, we recall that, by virtue of Proposition 1.2 in
Chap. 6 of \cite\zaktwo, a non-singular projective variety $X$ in
$\Bobb P^m \Bobb C$ which does not lie in a hyperplane is a {\it
Scorza\/} variety if and only if it satisfies (4.3) and (4.4)
below:
\newline\noindent
(4.3) $k_0(X) = \left[ \frac {\dim X} \delta \right]$;
\newline\noindent
(4.4) $\delta_i = i \delta$ for $1 \leq i \leq k_0(X)$.

With these preparations out of the way, we recall Zak's
classification, cf. Theorem 5.6 in Chap. 6 of \cite\zaktwo: For
the case where $\delta$ divides $\dim X$, given $\sk \geq 2$, the
$\sk$-Scorza varieties $X$ are exactly the projectivizations of
rank 1 strata in the simple complex Jordan algebras of rank $k+1$,
that is, the regular ones listed above as (1.1.$\sk$),
(1.2.$\sk$.r), (1.3.$\sk$.r), together with the Severi variety
(1.4) in case $\sk=2$. The classification of regular Scorza
varieties has been reworked in \cite\chapuone. When $\delta$ does
not divide $\dim X$, for $\sk \geq 2$, there remain two families
of $\sk$-Scorza varieties $X$ which are exactly those listed above
as (1.2.$\sk$.n) and (1.3.$\sk$.n); these  are the
projectivizations of rank 1 strata in the corresponding rank $k+1$
hermitian Jordan triple systems $\roman M_{\sk+1,\sk+2}(\Bobb C)$
and $\Lambda^2(\Bobb C^{2\sk +3})$, respectively. In particular,
given an $n$-dimensional non-singular projective variety $X$,
since $k_0(X) \leq \frac n \delta$, the requirement $k_0(X) =
\left[ \frac n \delta \right]$ is equivalent to $S^{\left[ \frac n
\delta \right]-1}X \ne \Bobb P^m \Bobb C$ and, when $X$ is a {\it
regular\/} $\sk$-Scorza variety ($\sk \geq 2$), $S^{\sk-1}X$ is a
hypersurface of degree $\sk+1$. More details about this
hypersurface will be given in Remark 5.6 below.

\medskip\noindent
{\bf 5. Holomorphic nilpotent orbits and the proof of Theorem 1.5}
\smallskip\noindent
A detailed account of holomorphic nilpotent orbits is given
in our paper \cite\kaehler.
Here we recall only what we need for ease of exposition.

Let $(\fra g, z)$ be a semisimple Lie algebra of hermitian type,
with Cartan decomposition $\fra g = \fra k \oplus \fra p$.
We refer to an  adjoint orbit $\Cal O \subseteq \fra g$
having the property that the projection map from $\fra g$ to $\fra p$,
restricted to $\Cal O$, is a diffeomorphism onto its image,
as a {\it pseudoholomorphic\/} orbit. A pseudoholomorphic orbit $\Cal O$
inherits a complex structure from the complex structure $J_z$ on $\fra p$,
and this complex structure, combined with the Kostant-Kirillov-Souriau
form on $\Cal O$, viewed as a coadjoint orbit by means of
(a positive multiple of) the Killing form, turns
$\Cal O$ into a (not necessarily positive) K\"ahler manifold.

We now choose a positive multiple of the Killing form. We say that
a pseudoholomorphic orbit $\Cal O$ is {\it holomorphic\/} provided
the resulting K\"ahler structure on $\Cal O$ is positive. The name
\lq\lq holomorphic\rq\rq\ is intended to hint at the fact that
the holomorphic discrete series representations of $G$ arise from holomorphic
quantization on integral {\it semisimple\/} holomorphic orbits but, beware,
the requisite complex structure (needed for the construction of the
holomorphic discrete series representation) is not the one arising from
projection to $\fra p$.

Let $\Cal O$ be a holomorphic {\it nilpotent\/} orbit, and let
$C^{\infty}(\overline{\Cal O})$ be the algebra of {\it Whitney smooth
functions\/} on the ordinary topological closure $\overline{\Cal O}$ of
$\Cal O$ (not the Zariski closure) resulting from the embedding of
$\overline{\Cal O}$ into $\fra g^*$. The Lie bracket on $\fra g$ passes to
a Poisson bracket $\{\cdot,\cdot\}$ on $C^{\infty}(\overline{\Cal O})$,
even though $\overline{\Cal O}$ is {\it not\/} a smooth manifold. This
Poisson bracket turns $\overline{\Cal O}$ into a stratified symplectic space.

A {\it complex analytic stratified K\"ahler structure\/} on a
stratified space $N$ is a stratified symplectic structure
$(C^{\infty}(N),\{\cdot,\cdot\})$ ($N$ is not necessarily smooth
and $C^{\infty}(N)$ not necessarily an algebra of ordinary smooth
functions) together with a complex analytic structure which, on
each stratum, \lq\lq combines\rq\rq\  with the symplectic
structure on that stratum to a K\"ahler structure, in the
following sense: The stratification underlying the stratified
symplectic structure is a refinement of the complex analytic
stratification whence each stratum is a complex manifold; each
holomorphic function is a smooth function in $C^{\infty}(N,\Bobb
C)$; and on each stratum, the Poisson  structure is symplectic in
such a way that the symplectic structure combines with the complex
structure to a K\"ahler structure. See Section 2 of \cite\kaehler\
for details. The structure may be described in terms of a {\it
stratified K\"ahler polarization\/} \cite\kaehler\ which induces
the K\"ahler polarizations on the strata and encapsulates the
mutual positions of these polarizations on the strata. A complex
polarization can no longer be thought of as being given by the
$(0,1)$-vectors of a complex structure, though. When the complex
analytic structure is normal we simply refer to a {\it normal
K\"ahler structure\/}. We now recall a few facts from
\cite\kaehler.

\noindent
{(5.1)} {\sl Given a holomorphic nilpotent orbit $\Cal O$, the closure
$\overline{\Cal O}$ is a union of finitely many holomorphic nilpotent orbits.
Moreover, the diffeomorphism from $\Cal O$ onto its image in $\fra p$
extends to a homeomorphism from the closure $\overline{\Cal O}$
onto its image in $\fra p$, this homeomorphism turns $\overline{\Cal O}$
into a complex affine variety, and the complex analytic structure, in turn,
combines with the  Poisson structure
$(C^{\infty}(\overline{\Cal O}),\{\cdot,\cdot\})$
to a normal (complex analytic stratified) K\"ahler structure.\/}
See Theorem 3.2.1 in \cite\kaehler.

Let  $r$ be the real rank of $\fra g$. There are $r+1$ holomorphic nilpotent
orbits $\Cal O_0,\dots, \Cal O_r$, and these are linearly ordered
in such a way that
$$
\{0\}=\Cal O_0 \subseteq \overline{\Cal O_1}
\subseteq \ldots \subseteq \overline{\Cal O_r} = \fra p,
\tag5.2
$$
cf. (3.3.10) in \cite\kaehler. Recall that the Cartan decomposition induces the
decomposition
$\fra g^{\Bobb C} = \fra k^{\Bobb C}\oplus \fra p^+\oplus \fra p^-$
of the complexification $\fra g^{\Bobb C}$ of $\fra g$,
$\fra p^+$ and $\fra p^-$ being the holomorphic and antiholomorphic
constituents, respectively, of $\fra p^{\Bobb C}$.

\noindent
{(5.3)} {\sl The projection from $\overline {\Cal O_r}$ to $\fra p$ is a
homeomorphism onto $\fra p$, and the $G$-orbit stratification of
$\overline {\Cal O_r}$ passes to the $K^{\Bobb C}$-orbit stratification
of $\fra p\cong \fra p^+$. Thus, for $1 \leq s \leq r$, restricted to
$\Cal O_s$, this homeomorphism is a $K$-equivariant diffeomorphism
from $\Cal O_s$ onto its image in $\fra p^+$, and this image is a
$K^{\Bobb C}$-orbit in $\fra p^+$.\/} See Theorem 3.3.11 in \cite\kaehler.

\noindent {\smc Remark 5.4.} The holomorphic nilpotent orbits in a
simple hermitian Lie algebra $\fra g$ are precisely those which
have the property that the projection to $\fra p\cong\fra p^+$
realizes the Kostant-Sekiguchi correspondence. See Remark 3.3.13
in \cite\kaehler.

\demo{Proof of Theorem {\rm 1.5}} The ascending sequence (5.2)  of
affine complex varieties determines the ascending sequence
$$
Q_1 \subseteq Q_2 \subseteq  \dots \subseteq Q_r=\Bobb P(\fra p)
\tag5.5
$$
of projective varieties where $\Bobb P(\fra p)$ is the projective
space on $\fra p$ and where each $Q_s$ ($1 \leq s \leq r$) arises
from $\overline{\Cal O_s}$ by projectivization. The stratified
K\"ahler structures on the constituents of (5.2) pass to
stratified K\"ahler structures on the constituents of (5.5), and
all stratified K\"ahler structures in sight are normal. See
Section 10 in \cite\kaehler\ for details.

Exploiting the sequence (5.5) for each of the four simple regular
rank 3 hermitian Lie algebras, from  the description of Severi
varieties in terms of the simple regular rank 3  hermitian Lie
algebras given in Section 4, we conclude that Theorem 1.5 holds
for $k=2$ in the four regular cases or, equivalently, that the
statements in the Addendum 1.7 hold. In particular, the discussion
in Section 8 (Appendix) of \cite\atiybern\ confirms that, under
these circumstances, $Q_2$ is indeed the chordal variety of the
Severi variety $Q_1$, the exceptional case (1.4) being included.
Likewise, for arbitrary rank $r\geq 3$, exploiting the sequence
(5.5) for each of the simple rank $r$ hermitian Lie algebras, from
the description of $(r-1)$-Scorza varieties in terms of the simple
regular  hermitian Lie algebras of rank $r$ given in Section 4, we
deduce that Theorem 1.5 holds for arbitrary rank $r\geq 3$ as
well, the non-regular cases being included. In particular, in the
remaining cases, which are all classical, that is, the
corresponding Jordan algebras or JTS's are matrix Jordan algebras
or matrix JTS's over the associative division algebras, any of the
corresponding Scorza varieties is the corresponding subspace $Q_1$
arising from Jordan rank 1 matrices by projectivization and, for
$1 \leq k\leq r-1$, the secant variety $S^k(Q_1)$ arises likewise
from the closure of the corresponding stratum of Jordan rank
$k+1$; indeed, this stratum consists of ordinary matrices of
Jordan rank $k+1$, and this observation entails the secant
properties. \qed
\enddemo

\noindent {\smc Remark 5.6.} In each of the four regular rank 3
cases (critical Severi varieties), the cubic polynomial which
determines the hypersurface $Q_2$ is the fundamental relative
invariant (Bernstein-Sato polynomial) of the corresponding
(irreducible regular) pre-homogeneous space, cf. Theorem 7.1 in
\cite\kaehler\ for the cases (1.1)--(1.3) and Theorem 8.4.1 in
\cite\kaehler\ for the case (1.4); in the latter case, the cubic
polynomial is the generic norm or Freudenthal's generalized
determinant mentioned in Section 4. Likewise, still in view of
Theorem 7.1 in \cite\kaehler, for higher rank $r \geq 4$, for each
of the three regular $(r-1)$-Scorza varieties $Q_1$, the degree
$r$ polynomial which determines the hypersurface $Q_{r-1} =
S^{r-2}Q_1 $ is the fundamental relative invariant (Bernstein-Sato
polynomial) of the corresponding (irreducible regular)
pre-homogeneous space.

\noindent {\smc Remark 5.7.} Let $p\geq q\geq 3$. The hermitian
Jordan triple system $\roman M_{q,p}(\Bobb C)$ has rank $q$ and
gives rise to a sequence
$$
\Bobb P^{q-1}\Bobb C \times \Bobb P^{p-1}\Bobb C = Q_1 \subseteq
Q_2 \subseteq  \ldots \subseteq Q_q=\Bobb P^{pq-1}\Bobb C
\tag5.7.1
$$ of stratified K\"ahler spaces of the kind (1.6) which
includes the Segre embedding of $Q_1=\Bobb P^{q-1}\Bobb C \times
\Bobb P^{p-1}\Bobb C $ into $\Bobb P^{pq-1}\Bobb C$. In this case,
$\delta =2$,
$$
k_0(Q_1)+ \left[\frac {p-q}2\right] = \left[ \frac {\dim Q_1}
\delta \right],
$$
and the requirement (4.4) holds (with $Q_1$ substituted for $X$),
whence $Q_1$ is a $(q-1)$-Scorza variety if and only if $p=q$ or
$p=q+1$; indeed, these varieties are then exactly those which
belong to the family (1.2.$\sk$.n). Even though, for $p \geq q+2$,
$\Bobb P^{q-1}\Bobb C \times \Bobb P^{p-1}\Bobb C$ is not a Scorza
variety, the assertions of Theorem 1.5 still hold for this case
(except the statement referring explicitly to a Scorza variety),
and the sequence (5.7.1) still enjoys similar secant properties as
the sequence (1.6) for an ordinary Scorza variety. The
classification of Scorza varieties together with this additional
class of varieties is {\it exactly parallel to the classification
of positive definite hermitian\/} JTS's and hence to the
classification of irreducible hermitian symmetric spaces of
non-compact type (which goes back to E. Cartan), cf. e.~g.
\cite\satakboo, and therefore as well to the classification of
irreducible bounded symmetric domains.

\noindent {\smc Remark 5.8.} By the same procedure as for higher
rank, for $q \geq 3$, the regular rank 2 simple hermitian Lie
algebra $\fra {so}(2,q)$ leads to the standard quadric in $\Bobb
P^{q-1}\Bobb C$ (a conic when $q=3$), cf. Sections 6 and 10 of
\cite\kaehler. Since $\fra {sp}(2,\Bobb R)\cong \fra {so}(2,3)$,
$\fra {su}(2,2)\cong \fra {so}(2,4)$, $\fra{so}^*(8)\cong \fra
{so}(2,6)$, every regular simple rank 2 hermitian Lie algebra
occurs here. Likewise, the non-regular rank 2 simple hermitian Lie
algebras $\fra {su}(p,2)$ ($p \geq 3$), $\fra{so}^*(10)$,
 and $\fra e_{6(-14)}$ lead
to the projective varieties, respectively, $\Bobb P^{p-1}\Bobb C
\times \Bobb P^1\Bobb C \subseteq \Bobb P^{2p-1} \Bobb C$ (Segre
embedding), $G_2\Bobb C^5 \subseteq \Bobb P^{9} \Bobb C$
(Pl\"ucker embedding),
 and $Q \subseteq \Bobb P^{15}\Bobb C$,
where $Q$ arises from the subspace of pure spinors in $\fra p^+
\cong \Cal D^5_+$ (positive half-spin representation of complex
dimension 16) by projectivization; see Sections 8 and 10 in
\cite\kaehler. These varieties can be interpreted as the limiting
case of {\rm 1}-Scorza varieties; in particular, with $k=1$, the
assertions in Theorem 1.5 still hold for these varieties and, in
the regular case, these varieties are quadratic.

\medskip\noindent {\bf 6. The classical cases via K\"ahler reduction}
\smallskip\noindent
Given a Hodge manifold $N$, endowed with an appropriate group of
symmetries and momentum mapping, reduction carries it to a complex
analytic stratified K\"ahler space $N^{\roman{red}}$ which is as
well a projective variety \cite\kaehler\ (Section 4) and the
question arises whether a {\sl complex projective space into which
$N^{\roman{red}}$ embeds carries an exotic structure which, via
the (Kodaira) embedding, restricts to the complex analytic
stratified K\"ahler structure on\/} $N^{\roman{red}}$. We will now
show that, for the classical cases, each constituent $Q_s$ of the
ascending sequence (1.6) $(1 \leq s \leq r =k+1\geq 3)$ is of this
kind. The construction to be given below includes, in particular,
Zak's Veronese mapping \cite\zaktwo, cf. Remark 6.8 below.

Let $s \geq 1$, $\din \geq 1$, let $\Bobb K = \Bobb R, \Bobb C,
\Bobb H$, and consider the standard (right) $\Bobb K$-vector space
$\Bobb K^s$ of dimension $s$, endowed with a (non-degenerate)
positive definite hermitian form $(\cdot,\cdot)$; further, let $V
= \Bobb K^{\din}$,  endowed with a skew form $\Cal B$ and
compatible complex structure $J_V$  such that associating $\Cal
B(u,J_V v)$ to $u,v \in V$ yields a positive definite hermitian
form on $V$. Moreover, let $\GHH= \GHH(s)= \roman
U(V^s,(\cdot,\cdot))$, $\fra h = \roman{Lie}(\GHH)$, $\HG = \roman
U(V,\Cal B)$, $\fra g = \roman{Lie}(\HG)$, and denote the split
rank of $\HG = \roman U(V,\Cal B)$ by $r$. More explicitly:

\noindent (6.1) $\Bobb K=\Bobb R,\ V = \Bobb R^{\din},\ {\din} = 2
\ell,\ \GHH= \roman O(s),\ \HG= \roman {Sp}(\ell,\Bobb R),\ r =
\ell$;
\newline\noindent
(6.2) $\Bobb K=\Bobb C,\ V = \Bobb C^{\din},\ {\din} = p +q,\
\GHH= \roman U(s),\ \HG= \roman U(p,q)$, $p \geq q$, $r = q$;
\newline\noindent
(6.3) $\Bobb K=\Bobb H,\ V = \Bobb H^{\din},\ \GHH= \roman
U(s,\Bobb H) = \roman {Sp}(s),\ \HG= \roman O^*(2 {\din})$, $r =
[\frac {\din}2]$.

Let $W=W(s)=\roman{Hom}_{\Bobb K}(\Bobb K^s,V)$. The 2-forms
$(\cdot,\cdot)$ and $\Cal B$ induce a symplectic structure $\omega_W$ on $W$,
and the groups $\KK$ and $\HG$ act on $W$ in an obvious fashion:
Given $x \in \KK, \ \alpha \in \roman{Hom}_{\Bobb K}(\Bobb K^s,V),\ y \in \HG$,
the action is given by the assignment to $(x,y,\alpha)$ of $y \alpha x^{-1}$.
These actions preserve the symplectic structure $\omega_W$ and are hamiltonian:
Given $\alpha \colon \Bobb K^s \to V$, define
$\alpha^\dagger \colon V \to \Bobb K^s$ by
$$
(\alpha^\dagger \bold u,\bold v) = \Cal B(\bold u,\alpha \bold v),
\ \bold u \in V,\, \bold v \in \Bobb K^s.
$$
Under the identifications of $\fra {\kk}$ and $\fra {\hg}$ with their duals
by means of the half-trace pairings, the momentum mappings $\mu_{\KK}$ and
$\mu_{\HG}$ for the $\KK$- and $\HG$-actions on $W$ are given by
$$
\align
\mu_{\KK} \colon W &@>>> \fra {\kk},
\quad
\mu_{\KK}(\alpha)= -\alpha^\dagger \alpha \colon \Bobb K^s \to \Bobb K^s,
\\
\mu_{\HG} \colon W &@>>> \fra {\hg},
\quad
\mu_{\HG}(\alpha) =
 \alpha \alpha^\dagger \colon V \to V
\endalign
$$
respectively; these are the  Hilbert maps of invariant theory for the
$\KK$- and $\HG$-actions. Moreover, the groups $G$ and $H$ constitute
a reductive dual pair in $\roman{Sp}(W,\omega_W)$, $J_V$ lies in
$\fra {\hg}$, and $(\fra g,\frac 12 J_V)$ is a simple (or reductive with
simple semisimple constituent) Lie algebra of hermitian type. For
$\phi \colon \Bobb K ^s \to V$, let $J_W (\phi) = J_V \circ \phi$.
This yields a complex structure $J_W$ on $W$ which,
together with the symplectic structure $\omega_W$, turns $W$ into a flat
(not necessarily positive) K\"ahler manifold.
When we wish to emphasize in notation that $W$ is endowed
with this complex structure we write $W_J$.
The $H$-action on $W_J$ preserves the complex structure.
See Section 5 of \cite\kaehler\ for details.

Since the $G$- and $H$-actions centralize each other, the $G$-momentum mapping
$\mu_{\HG}$ induces a map $\overline \mu_{\HG}$ from the ${\KKK}$-reduced space
$W(s)^{\roman{red}} = \mu_{\KKK}^{-1}(0)\big / {\KKK}$ into $\fra {\hg}$.
We now recall the following facts from Theorem 5.3.3 in \cite\kaehler.

\noindent (6.4)
{\rm (1)}
{\sl The induced map $\overline \mu_{\HG}$ from the ${\KKK}$-reduced space
$W(s)^{\roman{red}}$ into $\fra {\hg}$  is a proper embedding of
$W(s)^{\roman{red}}$ into $\fra {\hg}$ which, for $s \leq r$,
induces an isomorphism of normal K\"ahler spaces
from $W(s)^{\roman{red}}$ onto the closure $\overline {\Cal O_{s}}$
of the holomorphic nilpotent orbit\/} $\Cal O_{s}$.
\newline\noindent
{\rm (2)} {\sl For $2\leq s \leq r$, the injection of $W(s-1)$ into $W(s)$
induces an injection of $W(s-1)^{\roman{red}}$ into $W(s)^{\roman{red}}$
which, under the identifications with the closures of holomorphic
nilpotent orbits, amounts to the inclusion\/}
$\overline {\Cal O_{s-1}} \subseteq \overline {\Cal O_s}$.
\newline\noindent
{\rm (3)} {\sl For $s > r$,
the obvious injection of $W(s-1)$ into $W(s)$ induces an isomorphism
of $W(s-1)^{\roman{red}}$ onto $W(s)^{\roman{red}}$.\/}
\newline\noindent
{\rm (4)} {\sl For $s \geq 1$, under the projection from
$\fra g =\fra k \oplus \fra p$ to $\fra p$, followed by the identification
of the latter with $\fra p^+$, the image of the reduced space
$W^{\roman{red}}$ in $\fra g$ is identified with the affine complex
categorical quotient $W_J\big /\big / {\KKK}^{\Bobb C}$, realized in
the complex vector space $\fra p^+$.\/}

The projective space $\Bobb P[W]$ being endowed with a positive multiple of
the Fubini-Study metric, the momentum mapping passes to a momentum mapping
$\mu_{\KK} \colon \Bobb P[W] \to \fra {\kk}$. For $s=r$,
K\"ahler reduction yields the projective space $\Bobb P[\fra p]$ on $\fra p$,
endowed with an exotic K\"ahler structure;
likewise, for $1 \leq s \leq r$, K\"ahler reduction yields a
normal K\"ahler space $Q_s$ together with an ascending sequence
$$
Q_1 \subseteq Q_2 \subseteq \ldots\subseteq  Q_s
\tag6.5
$$
of normal K\"ahler spaces which are, in fact, the closures of the
strata of $Q_s$; the statement (6.4) above implies that this
sequence amounts to the sequence (5.5), truncated at $Q_s$.
Moreover, the embeddings of the $Q_s$'s into $Q_r=\Bobb P[\fra p]$
are Kodaira embeddings. In view of (6.4) above, the normal
K\"ahler structure coming from K\"ahler reduction coincides with
the structure coming from projectivization of the closure of the
corresponding holomorphic nilpotent orbit. See Section 10 of
\cite\kaehler\ for details.

Let $r\geq 3$ and $s=r$; with $\din =2r$, the sequence (6.5) then
recovers the ascending sequence (1.6) for each of the three
classical rank $r$ cases, where the parameter $n$ in (1.6) is now
written as $r$. In particular, for $r=3$, we obtain the sequence
(1.6.SEVERI) for the three regular classical cases. For
illustration, we describe briefly the various constituents in this
particular case where the numbering (6.6.*) corresponds to the
numbering (1.*) for $*=1,2,3$:
\newline\noindent
The sequence $\Bobb P[W(1)] \subseteq \Bobb P[W(2)] \subseteq \Bobb P[W]$ has
the form:
\newline\noindent
(6.6.1)
$\Bobb P^2 \Bobb C \subseteq \Bobb P^5 \Bobb C \subseteq \Bobb P^8 \Bobb C$
\newline\noindent
(6.6.2)
$\Bobb P^5 \Bobb C \subseteq \Bobb P^{11} \Bobb C \subseteq\Bobb P^{17}\Bobb C$
\newline\noindent
(6.6.3)
$\Bobb P^{11}\Bobb C\subseteq\Bobb P^{23}\Bobb C\subseteq\Bobb P^{35}\Bobb C$
\newline\noindent
The sequence (1.6.SEVERI) has the form:
\newline\noindent
(6.6.1$'$)
$X =\Bobb P^2\Bobb C \subseteq Q^4 \subseteq \Bobb P^5 \Bobb C
\cong\Bobb P^8 \Bobb C \big/\big/ \roman O(3,\Bobb C)$
\newline\noindent
(6.6.2$'$)
$X=\Bobb P^2\Bobb C\times\Bobb P^2\Bobb C\subseteq Q^7\subseteq\Bobb P^8\Bobb C
\cong\Bobb P^{17} \Bobb C \big/\big/ \roman {GL}(3,\Bobb C)$
\newline\noindent
(6.6.3$'$)
$X=\roman G_2(\Bobb C^6) \subseteq Q^{13}\subseteq\Bobb P^{14}\Bobb C
\cong\Bobb P^{35} \Bobb C \big/\big/ \roman {Sp}(3,\Bobb C)$
\newline\noindent
As GIT-quotients, the Severi varieties $X$ and the cubics $Q^*$
may be written out in the following fashion:
\newline\noindent
(6.6.1$''$)
$X=\Bobb P^2 \Bobb C \big/\big/ \roman {O}(1,\Bobb C) \cong \Bobb P^2\Bobb C$,
$Q^4 =\Bobb P^5 \Bobb C \big/\big/ \roman O(2,\Bobb C)$
\newline\noindent
(6.6.2$''$)
$X=\Bobb P^2\Bobb C\times\Bobb P^2\Bobb C
\cong
\Bobb P^5 \Bobb C \big/\big/ \roman {GL}(1,\Bobb C)$,
$Q^7 =
\Bobb P^{11} \Bobb C \big/\big/ \roman {GL}(2,\Bobb C)$,
\newline\noindent
(6.6.3$''$)
$X=\roman G_2(\Bobb C^6)
\cong
\Bobb P^{11} \Bobb C \big/\big/ \roman {Sp}(1,\Bobb C)$,
$Q^{13} =\Bobb P^{23} \Bobb C \big/\big/ \roman {Sp}(2,\Bobb C)$.
\newline\noindent
Here, for $1 \leq s \leq 3$, the actions of the groups $\roman
O(s,\Bobb C)$, $\roman {GL}(s,\Bobb C)$, and $\roman {Sp}(s,\Bobb
C)$ on the corresponding projective spaces arise from the actions
of the groups $H$ on $W(s)$ listed in (6.1)--(6.3) above. The
sequences (6.6.1)--(6.6.3) may be viewed as resolutions of
singularities (in the sense of stratified K\"ahler spaces) for the
corresponding sequences (1.6.SEVERI). The disjoint union
$$
\Bobb P[W]= H^{\Bobb C}\Bobb P[W(1)] \cup
H^{\Bobb C}(\Bobb P[W(2)] \setminus \Bobb P[W(1)])\cup
(\Bobb P[W] \setminus H^{\Bobb C}\Bobb P[W(2)])
$$
is the  $H^{\Bobb C}$-orbit type decomposition of $\Bobb P[W]$ in each case;
here $H^{\Bobb C}$ denotes the complexification of $H$.

In particular, the cubic $Q_2$
(written as $Q^4$, $Q^7$, $Q^{13}$ according to the case considered
where the superscript indicates the complex dimension)
arises by K\"ahler reduction, applied to the projective space $\Bobb P[W(2)]$.
This is an instance of the situation referred to at the beginning
of this section. Furthermore, the Severi variety $Q_1$
arises by K\"ahler reduction, applied to
the projective space $\Bobb P[W(1)]$
on  $W(1) (=V)$.

\noindent
{\smc Remark 6.7.}
In the special case where
$(\fra g, \fra h) = (\fra {sp}(\ell,\Bobb R), \fra {so}(s,\Bobb R))$,
the space $W$ may be viewed as the unreduced phase space of $\ell$
particles in $\Bobb R^s$, and $\mu_{\GH}$ is the angular momentum mapping;
cf. \cite\lermonsj.
Thus the normal K\"ahler space $\overline {\Cal O_s}$ arises as the reduced
phase space of a system of $\ell$ particles in $\Bobb R^s$ with total
angular momentum zero. Likewise, the compact normal K\"ahler space $Q_s$
which, complex analytically, is a projective variety, arises as the reduced
phase space of a system of $\ell$ harmonic oscillators in $\Bobb R^s$ with
total angular momentum zero and constant energy. Here the energy is encoded
in the stratified symplectic Poisson structure; changing the energy amounts
to rescaling the Poisson structure. The constituents $Q_s$
($1 \leq s \leq \ell$) of the ascending sequence (5.5) have the following
interpretation: The top stratum $Q_{\ell} \setminus Q_{\ell-1}$ consists of
configurations in general position, that is, $\ell$ harmonic oscillators in
$\Bobb R^{\ell}$ with total angular momentum zero
such that the positions and momenta
do not lie in a cotangent bundle $\roman T^* \Bobb R^s$ for some
$\Bobb R^s \subseteq \Bobb R^{\ell}$ with $s <\ell$. For $1 \leq s<\ell$,
the stratum $Q_s \setminus Q_{s-1}$ consists of configurations of $\ell$
harmonic oscillators in $\Bobb R^{\ell}$
with total angular momentum zero such that the
positions and momenta lie in the cotangent bundle $\roman T^* \Bobb R^s$
for some $\Bobb R^s \subseteq \Bobb R^{\ell}$ but {\it not\/} in a
cotangent bundle of the kind $\roman T^* \Bobb R^{s-1}$, whatever
$\Bobb R^{s-1} \subseteq \Bobb R^{\ell}$; here the convention is that
$Q_0$ is empty. Under these circumstances, the complex analytic structure
does not have an immediate mechanical significance but helps understanding
the geometry of the reduced phase space. The complex analytic structure is
important for issues related with quantization, cf. \cite\kaehredu.
Though this is not relevant here, we note that,
for $s>\ell$, the obvious map from $Q_{\ell}$ to $Q_s$ is an isomorphism of
stratified K\"ahler spaces, and no new geometrical phenomenon occurs.
Since the other dual pairs lie in some $\fra{sp}(n,\Bobb R)$,
it is likely that in the other cases (corresponding to (1.2) and (1.3)),
the constituents $Q_s$ ($1 \leq s \leq r$) of the ascending sequence
(5.5) admit as well interpretations in terms of suitable constrained
mechanical systems.

\noindent {\smc Remark 6.8.} For $s=1$, the composite of the
momentum mapping $\mu_G\colon W \to \fra g$ with the orthogonal
projection to the symmetric part $\fra p$ of the Cartan
decomposition $\fra g = \fra k \oplus \fra p$ is exactly the {\it
Veronese mapping\/} introduced in Zak's book \cite\zaktwo.

\noindent {\smc Remark 6.9.} The exceptional Severi variety (1.4)
is notably absent here. The question whether this variety and the
corresponding ambient spaces $Q_2$ and $Q_3$  arise from K\"ahler
reduction in the same way as the varieties (1.1)--(1.3) and the
corresponding ambient spaces or, more generally, the varieties
(1.1.$\sk$)--(1.3.$\sk$.n) and the constituents of the
corresponding sequences of the kind (1.6), is related with that of
existence of a dual pair $(G,H)$ with $G=E_{7(-25)}$ and $H$
compact but apparently there is no such dual pair. It would be
extremely interesting to develop an alternative construction which
yields the ascending sequence (1.6.SEVERI) for the exceptional
Severi variety by K\"ahler reduction applied to a suitable
K\"ahler manifold, perhaps more complicated than just complex
projective space.

\bigskip
%\vfill \eject

\centerline{\smc References}
\medskip
\widestnumber\key{999}

\ref \no \armcusgo
\by J. M. Arms,  R. Cushman, and M. J. Gotay
\paper  A universal reduction procedure for Hamiltonian group actions
\paperinfo in: The geometry of Hamiltonian systems, T. Ratiu, ed.
\jour MSRI Publ.
\vol 20
\pages 33--51
\yr 1991
\publ Springer Verlag
\publaddr Berlin $\cdot$ Heidelberg $\cdot$ New York $\cdot$ Tokyo
\endref

\ref \no \atiybern
\by M. Atiyah and J. Berndt
\paper Projective planes, Severi varieties and spheres
\finalinfo{\tt math.DG/0206135}
\endref
\ref \no \chapuone \by P.-E.~Chaput \paper Scorza varieties and
Jordan algebras \jour Indagationes Math. \vol 14 \yr 2003 \pages
169--182 \finalinfo{\tt math.AG/0208207}
\endref
\ref \no \farakora
\by J. Faraut and A. Koranyi
\book Analysis On Symmetric Cones
\bookinfo Oxford Mathematical Monographs
 No. 114
\publ Oxford University Press
\publaddr Oxford U.K.
\yr 1994
\endref

\ref \no \freudone
\by H. Freudenthal
\paper Beziehungen der $E_7$ und $E_8$ zur Oktavenebene
\jour Indagationes Math.
\vol 16
\yr 1954
\pages 218--230
\endref

\ref \no \hartsboo
\by  R. Hartshorne
\book Algebraic Geometry
\bookinfo Graduate texts in Mathematics
 No. 52
\publ Springer
\publaddr Berlin-G\"ottingen-Heidelberg
\yr 1977
\endref

\ref \no \howeone
\by R. Howe
\paper Remarks on classical invariant theory
\jour  Trans. Amer. Math. Soc.
\vol 313
\yr 1989
\pages  539--570
\endref

\ref \no \poiscoho
\by J. Huebschmann
\paper Poisson cohomology and quantization
\jour
J. f\"ur die reine und angewandte Mathematik
\vol  408
\yr 1990
\pages 57--113
\endref
\ref \no  \souriau
\by J. Huebschmann
\paper On the quantization of Poisson algebras
\paperinfo Symplectic Geometry and Mathematical Physics,
Actes du colloque en l'honneur de Jean-Marie Souriau,
P. Donato, C. Duval, J. Elhadad, G.M. Tuynman, eds.;
Progress in Mathematics, Vol. 99
\publ Birkh\"auser Verlag
\publaddr Boston $\cdot$ Basel $\cdot$ Berlin
\yr 1991
\pages 204--233
\endref

\ref \no  \oberwork
\by J. Huebschmann
\paper Singularities and Poisson geometry of certain representation spaces
\paperinfo in: Quantization of Singular Symplectic Quotients,
N. P. Landsman, M. Pflaum, M. Schlichenmaier, eds.,
Workshop, Oberwolfach,
August 1999,
Progress in Mathematics, Vol. 198
\publ Birkh\"auser Verlag
\publaddr Boston $\cdot$ Basel $\cdot$ Berlin
\yr 2001
\pages 119--135
\finalinfo{\tt math.DG/0012184}
\endref

\ref \no \kaehler
\by J. Huebschmann
\paper K\"ahler spaces, nilpotent orbits, and singular reduction
\jour Memoirs AMS (to appear)
\finalinfo {\tt math.DG/0104213}
\endref

\ref \no \kaehredu \by J. Huebschmann \paper K\"ahler quantization
and reduction \jour J. f\"ur die reine und angewandte Mathematik
(to appear) \finalinfo {\tt math.SG/0207166}
\endref

\ref \no \lradq
\by J. Huebschmann
\paper Lie-Rinehart algebras, descent, and quantization
\jour Fields Institute Communications (to appear)
\finalinfo {\tt math.SG/0303016}
\endref

\ref \no \jacobstw
\by N. Jacobson
\paper Some projective varieties defined by Jordan algebras
\jour J. of Algebra
\vol 97
\yr 1985
\pages 556--598
\endref

\ref \no \kostasix
\by B. Kostant
\paper The principal three-dimensional subgroup
and the Betti numbers of a complex simple Lie group
\jour Amer. J. of Math.
\vol 81
\yr 1959
\pages 973--1032
\endref

\ref \no \lanmaone
\by J. M. Landsberg and L. Manivel
\paper The projective geometry of Freudenthal's magic square
\jour J. of Algebra
\vol 239
\yr 2001
\pages 477--512
\endref

\ref \no \lermonsj
\by E. Lerman, R. Montgomery, and R. Sjamaar
\paper Examples of singular reduction
\paperinfo Symplectic Geometry,
Warwick, 1990,  D. A. Salamon, editor,
London Math. Soc. Lecture Note
Series, vol. 192
\yr 1993
\pages  127--155
\publ Cambridge University Press
\publaddr Cambridge, UK
\endref
\ref \no \mccrione \by  K. McCrimmon \paper Jordan algebras and
their applications \jour Bull. Amer. Math. Soc. \vol 84 \yr 1978
\pages 612--627
\endref
\ref \no \mezzthom \by  E. Mezzeti and O. Tomassi \paper Some
remarks on varieties with degenerate Gauss image\jour Pac. J. of
Math. \vol 213\yr 2004 \pages 79--88
\endref

\ref \no \satakboo
\by I. Satake
\book Algebraic structures of symmetric domains
\bookinfo Publications of the Math. Soc. of Japan, vol. 14
\publ Princeton University Press
\publaddr Princeton, NJ
\yr 1980
\endref

\ref \no \scorzone \by G. Scorza \paper Determinazione delle
variet\`a  a tre dimensioni di $S_r$ $(r\geq 7)$ i ciu $S_3$
tangenti si tagliano a due a due \jour Rendiconti del Circolo
Matematico di Palermo \vol 25 \yr 1908 \pages 193--204
\finalinfo{Opera Soelte, vol. I, ed. Cremonese, Roma, 1960}
\endref

\ref \no \scorztwo \by G. Scorza \paper Sulle variet\`a a quattro
dimensioni di $S_r$ $(r\geq 9)$ i ciu $S_3$ tangenti si tagliano a
due a due \jour Rendiconti del Circolo Matematico di Palermo \vol
27 \yr 1909 \pages 148--178 \finalinfo{Opera Soelte, vol. I, ed.
Cremonese, Roma, 1960}
\endref

\ref \no \sekiguch
\by J. Sekiguchi
\paper Remarks on real nilpotent orbits of a symmetric pair
\jour J. Math. Soc. Japan
\vol 39
\yr 1987
\pages 127--138
\endref

\ref \no \severone
\by F. Severi
\paper Intorni ai punti doppi impropi di una superficie generale
dello spazio a quattro dimenzioni, e a'suoi punti tripli apparenti
\jour Rendiconti del Circolo Matematico di Palermo
\vol 15
\yr 1901
\pages 33--51
\endref

\ref \no \zakone
\by F. L. Zak
\paper Severi varieties
\jour Math. USSR, Sb.
\vol 54
\yr 1986
\pages 113--127
\endref

\ref \no \zaktwo
\by F. L. Zak
\book Tangents and Secants of Algebraic Varieties
\publ American Math. Soc.
\publaddr Providence, R. I.
\yr 1993
\endref

\enddocument